\documentclass[11pt,reqno,twoside]{article}
\usepackage{amsmath,amssymb,theorem}
\theorembodyfont{\itshape}
\newtheorem{thm}{Theorem}[section]
\newtheorem{lem}[thm]{Lemma}
\newtheorem{prop}[thm]{Proposition}
{\theorembodyfont{\upshape}
\newtheorem{define}[thm]{Definition}
\newtheorem{rem}[thm]{Remark}
\newtheorem{ex}[thm]{Example}
\newtheorem{ass}{Hypothesis}}
\newtheorem{cor}[thm]{Corollary}

\numberwithin{equation}{section}

\newcommand{\Proof}[1][]{\noindent{\itshape Proof#1: }}
\newcommand{\EndProof}{~$\Diamond$\bigskip}

\def\hat{\widehat}
\def\tilde{\widetilde}
\def\diam{\mathrm{diam\,}}
\def\SP{E}

\def\N{{\mathbb N}}
\def\R{{\mathbb R}}
\def\C{{\mathbb C}}
\def\E{{\mathbb E}}

\def\cB{\mathcal{B}}

\def\cF{\mathcal{F}}
\def\cG{\mathcal{G}}

\def\cN{\mathcal{N}}

\def\cX{\mathcal{X}}

\def\ua{\uparrow}

\def\ol{\overline}

\def\sm{\setminus}
\def \sub{\subset}
\def\Tr{{\rm Tr\,}}
\def\d{\delta}    
\def\a{\alpha}    \def\b{\beta} \def\g{\gamma}            \def\d{\delta}
\def\D{\Delta}     \def\ph{\varphi}       
  \def\g{\gamma}      \def\G{\Gamma}     
\def\l{\lambda}   \def\L{\Lambda}     \def\m{\mu}        \def\n{\nu}
\def\r{\rho}               
\def\p{\pi}              \def\s{\sigma}     \def\S{\Sigma}
          
\def\x{\xi}       
\def\z{\zeta}           

\newcommand{\CM}[1][\mu]{C_{#1}}
\newcommand{\CCM}[1][\mu]{\hat{C}_{#1}}

\begin{document}
\title{Conditional Intensity and Gibbsianness\\ of Determinantal
Point Processes}
\author{ Hans-Otto Georgii\footnotemark[1] \  and Hyun Jae
Yoo\footnotemark[2] }
\date{}
  \maketitle

\begin{abstract}
The Papangelou intensities of determinantal (or fermion) point
processes are investigated. These exhibit a monotonicity property
expressing the repulsive nature of the interaction, and satisfy a
bound implying stochastic domination by a Poisson point process.
We also show that determinantal point processes satisfy the
so-called condition $(\S_{\l})$ which is a general form of
Gibbsianness. Under a continuity assumption, the Gibbsian
conditional probabilities can be identified explicitly.
\end{abstract}
{\bf Running head}. {Determinantal point processes}

\footnotetext[1]{Mathematisches Institut der Universit\"{a}t
         M\"{u}nchen, Theresienstra{\ss}e 39, 80333 M\"{u}nchen, Germany.
         E-mail: georgii@lmu.de}
\footnotetext[2]{University College , Yonsei University, 134
Shinchon-dong, Seodaemoon-gu, Seoul 120-749, Korea. E-mail:
yoohj@yonsei.ac.kr}

\section{Introduction}

The aim of this paper is to investigate the dependence structure
of determinantal (also called fermion) point processes,
abbreviated DPP's. These are point processes on a locally compact
metric space $\SP$ with a particular repulsive dependence
structure induced by the fact that their correlation functions are
given by suitable determinants. More explicitly, the correlation
function $\r(\a)$ at a finite configuration $\a$ is the
determinant of the matrix obtained by evaluating a positive
definite integral kernel at the points of $\a$; see Subsection
\ref{subsec:DPP} for details.

Since their invention by Macchi \cite{M}, DPP's have attracted
much interest from various viewpoints. Spohn \cite{Sp1, Sp2}
investigated the dynamics of the so-called Dyson's model, a model
of interacting Brownian particles on the real line with pair force
$1/x$, or pair interaction $-\log x$. Its invariant measure is a
typical DPP having the sine kernel as defining integral kernel.
Rather recently, the theory of DPP's has been developed
much further. Soshnikov \cite{So} established the full existence
theorem for DPP's, discussed many examples occurring in various
fields of physics and mathematics,  and derived various further
results for DPP's with translation-invariant kernels, including
mixing properties and central limit theorems.
Shirai and
Takahashi \cite{ST1} also dealt with the existence theorem and
extended the theory to some generalized DPP's including
 Boson processes. They also established some particular properties
such as limit theorems, ergodic properties, and the Gibbs property
of the corresponding discrete model \cite{ST2}. In a series of
papers, Borodin and Olshanski studied the DPP's appearing in the
representation of the infinite symmetric group, see \cite{BO} and
the references therein. Lyons and Steif \cite{L, LS} investigated
the ergodic and stochastic domination properties of DPP's on
discrete lattices. The Glauber dynamics leaving some DPP's
invariant was studied by Shirai and the second author \cite{SY}.

In this paper we ask for the dependence properties of DPP's. Our leading
questions are the following:
\begin{itemize}
\item[--]What can be said about the repulsive nature of the interaction?
\item[--]When are determinantal point processes Gibbsian?
\end{itemize}
Unfortunately,  to answer these questions we need to exclude the
interesting case when $1$ belongs to the spectrum of the
underlying integral operator, and in general we can only establish
a weak form of Gibbsianness. To prove full Gibbsianness we need to
impose a natural though restrictive continuity assumption. A key
quantity we consider is the conditional intensity in the sense of
Papangelou,  which is a function $c(x,\x)$ of points $x\in\SP$ and
configurations $\x$. Its intuitive meaning is that, for a suitable
reference measure $\l$ on $\SP$,  $c(x,\x)\l(dx)$ is the
conditional probability of having a particle in $dx$ when the
configuration $\x$ is given. We will show that, locally on
relatively compact regions, Papangelou intensities always exist,
are given by ratios of determinants, and are decreasing functions
of $\x$ (Theorem~\ref{thm:PapMon}). This monotonicity is a natural
expression of the repulsiveness of DPP's. In particular, it
implies negative correlations of local vacuum events stating that
some local regions contain no particles
(Corollary~\ref{cor:Monoton}). We also show that the local
Papangelou intensities are bounded from above. This implies
stochastic domination of DPP's by Poisson point processes and, for
$\SP=\R^d$, the absence of percolation in the associated Boolean
model when the underlying integral kernel is small enough
(Corollaries~\ref{cor:PoissonDom} and \ref{cor:NoPerc}).

Next we will show that  the Papangelou
intensities of DPP's exist not only locally but also globally
(Theorem~\ref{thm:Sigma}). This means that all DPP's are Gibbsian
in a general sense, in that one can write down natural versions of
their conditional distributions in local regions when the
configuration outside of this region is fixed, and these
conditional distributions are absolutely continuous with respect
to the Poisson point process. It is a more difficult question to
decide whether the associated conditional Hamiltonians can be
expressed in terms of the underlying integral kernel in
the natural way one expects. Here we only have a partial
result, stating that this holds whenever these natural
versions are continuous almost everywhere (Theorem~\ref{thm:Gibbs}).
The last condition holds in particular when $\SP=\R^d$
and the associated interaction kernel has finite range and is
continuous and small enough (Proposition~\ref{prop:a.e.continuity}).

In the next section we set up the stage. That is, in
Subsection~\ref{subsec:DPP} we recall the definition of DPP's
together with some basic facts, while Subsection~\ref{subsec:PI}
contains a discussion of the Papangelou intensity and its
significance. Our results are stated in Section~\ref{sec:results}.
The proofs follow in Section~\ref{sec:proofs}. In the Appendix we
collect some useful facts and discuss in particular the question
of how to make a proper choice of integral kernels.

\vspace{0cm plus 4cm}

\section{Preliminaries}
\subsection{Determinantal point processes}
\label{subsec:DPP}

In this section we describe the general setting and recall the
definition of determinantal point processes (DPP). For a more
complete account of DPP's we refer to the survey \cite{So} and the
articles \cite{L, LS, ST1, ST2}.
 Let $E$ be any locally compact metric space serving as the state
space of the points, $\cB$ the Borel $\s$-algebra on $\SP$, and
$\cB_0$ the system of all \emph{relatively compact} Borel sets in
$\SP$. Also, let $\l$ be a diffuse, locally finite reference
measure on $(\SP,\cB)$. The standard case is when $\SP=\R^d$ and
$\l$ is Lebesgue measure.

Let $\cX$ be the space of all locally finite subsets (configurations) in
$\SP$, i.e.,
\[
\cX:=\{\x  \sub \SP:| \x \cap \L|<\infty \text{ for all
}\L \in \cB_0\}\,,
\]
where $|A|$ denotes the cardinality of a set $A$. Later on, we
will exploit the fact that $\cX$ is partially ordered by
inclusion. Given any $\L \sub \SP$, we write $\cX_\L:=\{\x \in\cX:
\x \sub\L\}$ for the set of all configurations in $\L$, and
$r_\L:\x \to\x_\L:= \x \cap\L$ for the corresponding projection
from $\cX$ onto $\cX_\L$. Also let $N_\L:\x \to | \x \cap \L|$ be
the associated counting variable on $\cX$, and $\cF_\L$ be the
smallest $\s$-algebra on $\cX$ such that $N_\D$ is measurable for
all relatively compact Borel sets $\D \sub \L$. We write $\cF$ for
$\cF_{\SP}$. Each configuration $\x \in \cX$ can be identified
with the integer-valued Radon measure $\sum_{x\in  \x }\d_x$ on
the Borel $\s$-algebra on $\SP$. The vague topology for the latter
then induces a topology on $\cX$ turning $\cX$ into a Polish
space. It is well-known that $\cF$ is the associated Borel
$\s$-algebra \cite{DV, Sh}. A \emph{point process} is a
probability measure $\m$ on $(\cX,\cF)$. We write $\m_\L:= \m\circ
r_\L^{-1}$ for its marginal on $\cX_\L$.

Next let $\cX_0=\{\a\in\cX: |\a|<\infty\}$ be the set of all
\emph{finite} configurations in $\SP$. $\cX_0$ is equipped with
the trace $\s$-algebra $\cF_0=\cF|_{ \cX_0}$ and the
\emph{$\l$-sample measure} $L$ defined by the identity
\begin{equation}\label{eq:L}
\int_{\cX_0} f(\a)\,L(d\a) =
\sum_{m=0}^{\infty}\frac{1}{m!}\int_{\SP^m} f(\{x_1,\ldots,
x_m\})\, \l(dx_1)\cdots \l(dx_m)
\end{equation}
for any measurable function $f:\cX_0\to\R_+$. For any $\L\subset\SP$ we
let $L_\L(d\a) = 1_{\{\a\subset\L\}}L(d\a)$ be the restriction of $L$ to
$\cX_\L$. Here  we use the notation $1_A$ for the indicator function of
a set $A$.

Recall that a point process $\m$ is said to have the
\emph{correlation function} $\r:\cX_0\to \R_+$ if $\r$  is
measurable and satisfies
\begin{equation}\label{eq:corrfct}
\int_\cX \  \sum_{ \a\in \cX_0:\, \a \sub \x} u( \a) \ \m(d \x )=
\int_{\cX_0}  u( \a) \r( \a)\,L(d \a)
\end{equation}
for any measurable $u:\cX_0\to \R_+$. (For a general account of
the interrelationship $\m\leftrightarrow \r$ between point
processes and correlation functions we refer to the recent paper
\cite{KK}.) The \emph{Poisson point process $\p^z$ with (locally
integrable) intensity function $z:\SP\to [0,\infty[$} is the
unique point process having correlation function $\r^z(\a
)=\prod_{x\in \a}z(x)$. Equivalently, $\p^z$ is the unique point
process such that, for each $\L \in\cB_0$, the projection
$\p^z_\L$ has the Radon-Nikodym density $\x\to e^{-\int_\L
z\,d\l}\prod_{x\in \x}z(x)$ relative to $L_\L$. The characteristic
feature of the determinantal point processes to be considered here
is that their correlation function is given by suitable
determinants. Given a function $K:\SP\times\SP \to\C$, we write
\begin{equation}\label{eq:Kmatrix}
K(\a,\a)=\big(K(x,y)\big)_{x,y\in\a}
\end{equation}
for the matrix obtained by evaluating $K$ at a finite configuration
$\a\in\cX_0$.

\begin{define}\label{def:DPP} Let $K(x,y)$, $x,y\in \SP$, be the
integral kernel of a positive\footnote{We use the terms
``positive operator'' and ``positive definite
matrix'' always in the weak sense, in that $0$ may belong to the
spectrum.} Hermitian integral operator $K$ on $L^2(\SP,\l)$. A
point process $\m$ with correlation function
\[
\r(\a)=\det K(\a,\a)
\]
is called a \emph{determinantal point process} (abbreviated DPP)
with \emph{correlation kernel} $K$.
\end{define}

The DPP's defined above are also known as \emph{fermion point
processes}; see \cite{ST1} for the boson case where the
determinant is replaced by the permanent, and more general
determinantal processes. In fact, there exist interesting examples
of DPP's for non-Hermitian kernels such as the discrete Bessel
kernel; see \cite{BOO}, for example. In this paper, however, we
confine ourselves to the Hermitian case. On the technical side, we
note that the matrices \eqref{eq:Kmatrix} involve the values of
the integral kernel $K(x,y)$ not only for $\l^{\otimes 2}$-almost
all $(x,y)\in\SP^2$ but also on the diagonal of $\SP^2$, which is
a $\l^{\otimes 2}$-nullset because $\l$ is diffuse. So one needs
to make a proper choice of $K(x,y)$ which is \emph{positive
definite a.e.}, in that $K(\a,\a)$ is positive definite for
$L$-almost all $\a\in \cX_0$ (implying that the correlation
function $\r$ is indeed nonnegative $L$-a.e.). In Lemma
\ref{lem:loc_trace_formula} will be explained how this can be done
.

Concerning the existence of DPP's, we state the following theorem
from \cite{So}, see also \cite{M, DV, ST1}. $I$ stands for the
identity operator, and the ordering $S\le T$ between operators
means that $T-S$ is a positive operator.
\begin{thm}[Macchi, Soshnikov]\label{thm:frpf}
A Hermitian locally trace class operator $K$ on $L^2(\SP,\l)$
defines a DPP $\m$ if and only if $0\le K\le I$, and in this case
$\m$ is unique.
\end{thm}
 Any DPP $\m$ is locally
absolutely continuous with respect to the $\l$-sample measure $L$
and admits explicit expressions for the local densities. To be
specific, for each  $\L \in \cB_0$ let $P_\L: L^2(\SP,\l)\to
L^2(\L,\l_\L)$ be the projection operator and $K_\L:=P_\L KP_\L$
the restriction of $K$ onto $L^2(\L,\l_\L)$. That is, $K_\L$
admits the kernel $K_\L(x,y)=1_\L(x)K(x,y)1_\L(y)$. Suppose that
$\m$ is the unique DPP for an operator $K$ as in Theorem
\ref{thm:frpf}. Then the density function of $\m_\L$ with respect
to $L_\L$, the so-called Janossy density \cite{DV} of $\m$, is
given by \cite{ST1, So}
\begin{equation}\label{eq:density_function}
\s_\L(\x)=\det(I-K_\L)\det J_{[\L]}(\x,\x), \quad \x\in\cX_\L;
\end{equation}
here  $J_{[\L]}:=K_\L(I-K_\L)^{-1}$, and the normalization
constant $\det(I-K_\L)$ is to be understood as a Fredholm
determinant \cite{Si}. A priori, this formula requires that $0\le
K_\L<I$. However, one can show that \eqref{eq:density_function}
still makes sense when $\|K_\L\|=1$; see \cite[p.~935]{So}. In view
of \eqref{eq:density_function}, we call $J_{[\L]}$ the local
interaction operator and its kernel the local interaction kernel
for $\L$. Again, Lemma \ref{lem:positive_kernel} allows to choose
the kernel $J_{[\L]}$ in such a way that all matrices
$J_{[\L]}(\x,\x)$ are positive definite. Finally, we note that the
correlation function $\r$ can be recovered from the local
densities $(\s_\L)$ by
\[
\r(\a)=\int_{\cX_\L}\s_{\L}(\a\x)\,L_\L(d\x)
\]
for $\a\in \cX_\L$, where $\a\x$ is shorthand for $\a\cup\x$.

\subsection{Papangelou intensities}
\label{subsec:PI}

Here we summarize some facts concerning the reduced (compound)
Campbell measure of a point process $\m$ as well as its Papangelou
intensity which describes the local dependence of particles.
\begin{define}\label{def:redCampb}
(a) The \emph{reduced} (or modified) \emph{Campbell measure} of a
point process $\m$ is the measure $\CM$ on the product space
$(\SP\times\cX, \cB\otimes\cF)$ defined by
\[
\CM(A) = \int \m(d\x ) \sum_{x\in \x } 1_A(x, \x \sm x)\,, \quad A\in
\cB\otimes\cF,
\]
where $\x \sm x:=\x \sm\{ x\}$.

(b) The \emph{reduced compound Campbell measure} of $\m$ is the measure $\CCM$
on the product space $(\cX_0\times\cX, \cF_0\otimes\cF)$ satisfying
\[
\CCM(B) = \int \m(d\x ) \sum_{\a\in\cX_0:\,\a\subset \x } 1_B(\a, \x \sm \a)\,,
\quad B\in \cF_0\otimes\cF.
\]
\end{define}
It is well-known and easy to check that the reduced Campbell
measure of the Poisson point process $\p^z$ is given by
$\CM[\p^z](dx,d\x)=z(x) \l(dx)\p^z(d\x)$, and a classical result
of Mecke \cite{Mecke} states that $\p^z$ is the only point process
with this property. This fact suggests the following concept.

\begin{define}
A point process $\m$ is said to satisfy \emph{condition} $(\S_\l)$
if $\CM \ll \l\otimes\m$. Any Radon-Nikodym density $c$ of $\CM$
relative to $\l\otimes\m$ is then called (a version of) the
\emph{Papangelou (conditional) intensity} (abbreviated PI) of
$\m$.
\end{define}
More explicitly, $c$ is a PI of $\m$ if
\begin{equation}\label{eq:Papint}
\int \m(d\x ) \sum_{x\in \x } f(x, \x \sm x) = \int \l(dx) \int
\m(d\x )\, c(x,\x ) f(x,\x )
\end{equation}
for all measurable functions $f:\SP\times\cX \to \R_+$.
Another way of stating this is that the intensity
measure of $\m$ has a density $\r_1$ relative to $\l$, and the
\emph{reduced Palm distribution $\m^x$ of $\m$ at
$x$} (which is defined by the desintegration formula $\CM(dx,d\x) = \l(dx)\r_1(x)\,\m^x(d\x)$) is absolutely
continuous with respect to $\m$ with density $c(x,\cdot)/\r_1(x)$ for $\l$-a.a. $x\in\{\r_1>0\}$. Intuitively,
$c(x,\x )\l(dx)$ is the conditional probability for a particle in the differential region $dx$ given the
configuration $\x $. Also, if
\[
Gh(\x ) = \int \l(dx)\ c(x,\x ) \big[h(\x \cup x)-h(\x )\big] +
\sum_{x\in \x }\big[h(\x \sm x)-h(\x )\big]
\]
is the formal generator of a birth-and-death process with 
birth rate $c(x,\x )\l(dx)$ for a particle in $dx$ and
death rate $1$ for each particle, \eqref{eq:Papint} is 
equivalent to the reversibility equation $\int g\, Gh\,
d\m = \int h\, Gg\, d\m$; to see this it is sufficient to set 
$f(x,\x )=g(\x )h(\x \cup x)$.

The following remark lists some further consequences of condition $(\S_\l)$.
In particular, it shows that $(\S_\l)$ processes are Gibbsian in a
general sense.

\begin{rem}\label{rem:MWM}
For any $\m$ satisfying condition $(\S_\l)$ the following conclusions hold.

\smallskip
(a) The reduced compound Campbell measure $\CCM$ is absolutely
continuous with respect to $L\otimes\mu$ with a Radon-Nikodym density
$\hat c$ satisfying $\hat c(\emptyset,\xi) =1$ and
\begin{equation}\label{eq:compPapdens}
\hat c(\a,\x ) = c(x_1,\x) \prod_{i=2}^n
c(x_i,x_1 \ldots x_{i-1}\x) \quad\text{when }\a=\{x_1,\ldots, x_n\};
\end{equation}
here $x_1 \ldots x_{i-1}\x = \{x_1, \ldots, x_{i-1}\}\cup \x$. In
particular, the right-hand side of \eqref{eq:compPapdens} is almost
surely symmetric in $x_1,\ldots, x_n$. $\hat c$ is called the
\emph{compound Papangelou intensity} (CPI). Explicitly, the relation
$\hat c=d\CCM/d( L\otimes\mu)$   means that
\begin{equation}\label{eq:compPapint}
\int \m(d\x ) \sum_{\a\in\cX_0:\,\a\subset \x } f(\a, \x \sm \a) =
\int L(d\a) \int \m(d\x )\, \hat c(\a,\x )\, f(\a,\x )
\end{equation}
for any measurable $f:\cX_0\times\cX \to \R_+$, and follows easily from
\eqref{eq:Papint} by induction on $|\a|$.

\smallskip
(b) For an $f$ depending only on $\a$, a comparison of
\eqref{eq:compPapint} and \eqref{eq:corrfct} shows that the correlation
function and the CPI of $\m$ are related to each other by
\[
\r(\a) = \int \m(d\x)\, \hat c(\a,\x)\quad\text{for $L$-almost all $\a$.}
\]

\smallskip
(c) Let $\L\in\cB_0$. Applying \eqref{eq:compPapint} to $f(\a,\x)=
g(\a)h(\x) 1_{\{\a\subset\L, \,N_\L(\x)=0\}}$ for any
$\cF_\L$-measurable $g:\cX\to\R_+$ and $\cF_{\L^c}$-measurable
$h:\cX\to\R_+$ and taking conditional expectations we find
\[
\begin{split}
&\int \E_\m(g|\cF_{\L^c})\,h\, d\m = \int g h\, d\m
= \int f \,d\CCM\\
&=\int \m(d\x)\,h(\x)\  \m(N_\L{=}0\,|\,\cF_{\L^c})(\x)  \int L_\L(d\a)\,
g(\a)\,\hat c(\a,\x_{\L^c})
\end{split}
\]
Hence
\begin{equation}\label{eq:condexpect}
\E_\m(g|\cF_{\L^c})(\x) = \m(N_\L{=}0\,|\,\cF_{\L^c})(\x)\ \int g\; \hat
c(\cdot, \x_{\L^c}) \,dL_\L
\end{equation}
for $\m$-almost all $\x$ and any $\cF_\L$-measurable $g$. In particular,
for $g\equiv 1$ we find that $\m(N_\L{=}0\,|\,\cF_{\L^c})>0$ almost surely for
each $\L\in\cB_0$, a property introduced by Papangelou \cite{Papangelou}
as condition $(\S)$, and by Kozlov \cite{Kozlov} as the condition of
non-degenerate vacuum.
\end{rem}

The observations in the preceding remark are due to Matthes, Warmuth and
Mecke \cite[Section 3]{MWM} and give one part of their theorem below;
cf. also \cite[Theorem 2$'$\,]{NZ}.

\begin{thm}[Matthes, Warmuth and Mecke]\label{thm:MWM}
A point process $\m$ satisfies condition $(\S_\l)$ with PI $c$ if and
only if, for each $\L\in\cB_0$, $\m$ is absolutely continuous relative
to $L_\L\otimes\m_{\L^c}$ with a density  $\g_\L$ satisfying
$\g_\L(x\x)=0$ whenever $\g_\L(\x)=0$ and $x\in\L\sm\x$. In this case,
for $L_\L\otimes\m_{\L^c}$-almost all $\x \in\cX$ we have
\begin{equation}\label{eq:condPapdens}
\g_\L(\x ) = Z_\L(\x _{\L^c})^{-1}\; \hat c(\x_\L,\x _{\L^c})
\end{equation}
with $0< Z_\L(\x _{\L^c}) = \int \hat c(\cdot,\x _{\L^c})\,dL_\L
<\infty$, and $c(x,\x ) = \g_\L(x\x )/\g_\L(\x )$ for
$\l\otimes\m$-almost all $(x,\x )\in \L\times\cX$.
\end{thm}
More explicitly, Equation \eqref{eq:condPapdens} means that for
each bounded measurable function $f:\cX\to\R$, the conditional
expectation $\E_\m(f\,|\,\cF_{\L^c})(\x )$ has the version
\begin{equation}\label{eq:condProb}
G_\L(f\,|\,\x ) := Z_\L(\x _{\L^c})^{-1} \int f(\a\x _{\L^c})\, \hat
c(\a,\x _{\L^c})
L_\L(d\a).
\end{equation}

\begin{rem}\label{rem:locMWM}
Theorem \ref{thm:MWM} has a counterpart for the restriction of a
point process $\m$ to a local region $\L\in\cB_0$. Let $\l_\L$
be the restriction of $\l$ to the $\s$-algebra $ \cB_\L$ of Borel
subsets of $\L$. Then $\m_\L$ satisfies condition $(\S_{\l_\L})$
with a PI $c_\L$ if and only if $\mu_\L$ is absolutely continuous
with respect to $L_\L$ with a density $\s_\L$ having an increasing
zero-set $\{\s_\L=0\}$. In this case, $\s_\L =
\s_\L(\emptyset)\,\hat c_\L(\cdot,\emptyset)$ $L_\L$-almost
everywhere, and $c_\L(x,\xi) = \s_\L(x\x)/\s_\L(\x)$ for
$\l_\L\otimes\m_\L$ almost all $(x,\x)\in \L\times\cX_\L$. This
follows from the above by noting that $\m_\L$, regarded as a
measure on $(\cX, \cF)$ supported on $\{N_{\L^c}=0\}$, is trivial
on $\cF_{\L^c}$; cf. also \cite[Proposition 3.1]{GK}.
\end{rem}

\section{Results}
\label{sec:results}

Our results on DPP's are based on the following standing
assumption on the underlying correlation operator $K$.

\begin{ass}\label{hyp:H} $K$ is a Hermitian integral
operator on $L^2(\SP,\l)$ such that
\begin{enumerate}
\item $K$ is of local trace class, and \item $\text{spec\,}K\sub[0,1[\,$;
i.e., $0\le K\le I$ in the operator ordering, and $\|K\|<1$.
\end{enumerate}
We write $\m$ for the unique DPP with correlation kernel $K$.
\end{ass}

Let us comment on these assumptions. First,  $K$ is
necessarily a Carleman operator~\cite{HS}, in that its
integral kernel satisfies $K(x,\cdot)\in L^2(\SP,\l)$ for
$\l$-almost all $x\in E$; this will be shown in Lemma
\ref{lem:loc_trace_formula}(b). Secondly, it follows from (b)
that $J:=K(I-K)^{-1}$ exists as a bounded
operator; we call $J$ the \emph{(global) interaction operator}.
Since the Carleman operators form a right ideal
in the space of bounded operators on $L^2(\SP,\l)$
\cite[Theorem 11.6]{HS}, hypothesis (H) ensures that $J$ is
also a Carleman operator. Moreover, since $K$ is supposed
to be of local trace class, so is $J$; for, if $\L\in
\cB_0$ then $J_\L:=P_\L J P_\L\le (1-\|K\|)^{-1} K_\L$,
whence $\Tr J_\L\le (1-\|K\|)^{-1}\Tr K_\L<\infty$.
Assumption (a) implies further that, for each $\L\in \cB_0$,
the local operators $K_\L$ and $J_\L$ as well as
the local interaction operator $J_{[\L]}:=K_\L(I-K_\L)^{-1}$
satisfy the conditions of Lemma~\ref{lem:positive_kernel}.
So, their kernel functions can and will be chosen in such a
way that all associated matrices $K_\L(\a,\a)$, $J_\L(\a,\a)$
and $J_{[\L]}(\a,\a)$ defined in \eqref{eq:Kmatrix} are
positive definite. Finally, we note that necessarily
$0\in \text{spec\,}K$ because $K_\L$, being trace class, is
compact and thus has eigenvalues tending to zero.

Our first result describes the local behavior of DPP's under (H).
\begin{thm}\label{thm:PapMon}
For each $\L\in\cB_0$,
$\m_\L$ satisfies condition $(\S_{\l_\L})$, and a version
of its CPI $\,\hat c_\L$ is given by
\begin{equation}\label{eq:cLformula}
\hat c_\L(\a,\x)= \det J_{[\L]}(\a\x,\a\x)/\det J_{[\L]}(\x,\x)\,;
\end{equation}
the ratio is defined to be zero whenever the denominator vanishes.
This version also satisfies the inequalities
\begin{equation}\label{eq:cLinequality}
\hat c_\L(\a,\x)\ge\hat c_\L(\a,\eta)\quad\text{and}\quad 0\le \hat c_\L(\a,\x)\le
\det J_{[\L]}(\a,\a) \le \prod_{x\in\a}J_{[\L]}(x,x)
\end{equation}
whenever $\x\sub\eta$ and $\a\sub\L\sm\eta$.
\end{thm}
Here are some consequences of the theorem. Let us look first at
the local CPI $\,\hat c_\L$ in \eqref{eq:cLformula}.
Let $\L,\D\in\cB_0$ with $\L\subset\D$ be given,
and $f:\cX_\L\to\R_+$  an arbitrary
measurable function. Then, combining \eqref{eq:cLformula} with
Remark \ref{rem:locMWM} and \eqref{eq:condProb} we find that
\begin{equation}\label{eq:locCondProb}
\begin{split}
\E_\m(f\,|\,\cF_{\D\sm\L})(\x) &= G_{\L,\D}(f\,|\,\x)\\
&:= Z_{\L,\D}(\x_{\D\sm\L})^{-1} \int L_\L(d\a) \, f(\a)\,\hat
c_\D(\a,\x_{\D\sm\L})
\end{split}
\end{equation}
for $\m$-almost all $\x\in\cX$. Let us comment on these
conditional probabilities.

\begin{rem}\label{rem:G_L,D}
(a)  For any $\x\in\cX_{\D\sm\L}$, the normalization constant
$$
Z_{\L,\D}(\x)= \int L_\L(d\a) \,\hat
c_\D(\a,\x)
$$
is finite. In fact, $Z_{\L,\D}(\x)\le  \det(I+J_\L)$.
Moreover, $Z_{\L,\D}(\x)$ is non-zero whenever $\det
J_{[\D]}(\x,\x)>0$ because $L_\L(\{\emptyset\})=1$. This means
that $G_{\L,\D}(\,\cdot\,|\,\x)$ is a well-defined probability measure on $\cX_\L$
for \emph{all} such $\x$.

(b) For all $\x$ with $D:=\det J_{[\D]}(\x,\x)>0$, $G_{\L,\D}(\,\cdot\,|\,\x)$ 
is a DPP with interaction kernel
$J_{\L,\D}^\x(x,y):=D^{-1} \det J_{[\D]}(x\x,y\x)$, $x,y\in\L$. 
(The last matrix is defined in obvious analogy
to \eqref{eq:Kmatrix}.) This observation is analogous to 
the fact observed in \cite{ST1} that the reduced Palm
measures $\m^x$ of $\m$ are again determinantal with correlation kernel $K^x(y,y')=K(x,x)^{-1}\det K(yx,y'x)$.
\end{rem}

Turning to the inequalities \eqref{eq:cLinequality}, we emphasize that
the monotonicity of the functions $\hat c_\L(\a,\cdot)$ expresses the
repulsive nature of the particle interaction in an infinitesimal way.
As a matter of fact, this
monotonicity is analogous to the inequality
\[
\r(\a_1)\r(\a_2)\ge \r(\a_1\cup\a_2)\r(\a_1\cap\a_2),\quad
\a_1,\,\a_2\in \cX_0,
\]
derived in \cite{ST1} for the correlation function $\r$ of $\m$,
and follows in the same way from an inequality for determinants of
Hermitian matrices; see \eqref{eq:det3} below. The following
corollary provides an integral version of this repulsiveness.

\begin{cor}\label{cor:Monoton}
For any $\L\subset\D\in\cB_0$ and any measurable
function $f:\cX_\L\to\R_+$, the ratio $G_{\L,\D}(f\,|\,\x
)/G_{\L,\D}(N_\L{=}0\,|\,\x )$ is decreasing in $\x\in\cX_{\D\sm\L} $.
In particular, $G_{\L,\D}(N_\L{=}0\,|\,\cdot)$ is increasing, and for
$\L\sub\D_1\sub\D_2\in\cB_0$ we have
\begin{equation}\label{eq:Monoton}
\m(N_\L{=}0\,|\,N_{\D_2\sm\L}{=}0)\le \m(N_\L{=}0\,|\,N_{\D_1\sm\L}{=}0) \le
\m(N_\L{=}0)\,.
\end{equation}
\end{cor}

The statement of the corollary is weaker than one may hope. In
fact, one might guess that $G_{\L,\D}(f\,|\,\x )$ is decreasing in
$\x$ for any increasing $\cF_{\L}$-measurable function $f$. (The
corollary implies this assertion only for $f=1_{\{N_\L\ge 1\}}$).
However, we have some doubts whether this can be expected to hold
in general. For if $f$ depends only on some part $\L_0$ of $\L$,
then an increase of $\x$ may repel some particles from
$\L\setminus\L_0$, giving a chance to some additional particles in
$\L_0$, so that $G_{\L,\D}(f\,|\,\x )$ will increase. So, the
situation is less satisfactory for point processes than in the
discrete case; cf. Theorem 6.5 of \cite{L}. Nevertheless,
\eqref{eq:Monoton} implies that, for disjoint $\L,\D\in\cB_0$, the
events $\{N_\L=0\}$ and $\{N_\D=0\}$ are negatively correlated;
see Proposition 2.7 of \cite{ST2} for the corresponding result in
the discrete case.

Next we exploit the domination bound $c_\L(x,\x)\le
J_{[\L]}(x,x)$, which will imply that $\m$ is
stochastically dominated by a Poisson process.
Recall that a point process $\n$ is \emph{stochastically dominated}
by a point process $\n'$, written $\n\preceq\n'$, if
\[
\int f\,d\n\le \int f\,d\n'
\]
for every increasing measurable function $f:\cX\to\R$. Here, $f$
is said to be increasing if $f(\x)\le f(\eta)$ whenever
$\x\subset\eta$. The intensity function of the dominating Poisson
process shall be given by $z(x)=J(x,x)$, where $J(x,y)$ is a
\emph{properly chosen} integral kernel of $J$. To understand the
necessity of a proper choice one should note again that the
diagonal in $\SP^2$ is a $\l^{\otimes 2}$-nullset because $\l$ is
diffuse; the behavior of the integral kernel on the diagonal that
enters in the definition of $z(x)$ is therefore a priori unrelated
to the operator $J$ (except when the kernel is continuous which we
do not require here). However, a natural way of determining the
kernel function on the diagonal is in terms of the local trace
formula
\begin{equation}\label{eq:locTrace1}
\Tr\,(J_\L)=\int_\L J(x,x)\;\l(dx)\quad\text{for all }\L\in\cB_0,
\end{equation}
and such a choice is possible according to Lemma \ref{lem:loc_trace_formula}.
The function $z(x):=J(x,x)$  is then locally $\l$-integrable, and
the Poisson point process  $\p^{J}$ with intensity function $z(x)$ is well-defined.

\begin{cor}\label{cor:PoissonDom}
$\m\preceq\p^J$.
\end{cor}
In the case when $\SP=\R^d$ and $\l$ is Lebesgue measure, the last
result implies in particular that there is no percolation in the
Boolean model associated to $\m$ when $J(x,x)$ is small enough. Let
$b_R(x)$ denote the closed ball of radius $R<\infty$ in $\R^d$
centered at $x$, and for $\x \in\cX$ let
\[
B_R(\x )=\bigcup_{x\in \x }b_R(x)
\]
the associated Boolean set. $B_R(\x )$ splits into connected
components called \emph{clusters}. A cluster is called
\emph{infinite} if it consists of infinitely many points of $\x $,
or equivalently, if its diameter is infinite. It is well-known
\cite{G, MR} that, for $d\ge 2$, there exists a critical threshold
$0<z_c\equiv z_c(d,R)<\infty$ for Poisson percolation: For the
Poisson point process $\p^z$ with constant intensity function
$z>0$, one has
\begin{equation*}
\p^z\big(\,\exists \text{ infinite cluster of }B_R(\cdot)\,\big)=
\begin{cases} 0 &\text{if }z<z_c\,,\\1 &\text{if }z>z_c\,.\end{cases}
\end{equation*}
For $d=1$ one has $z_c=+\infty$ \cite[Theorem 3.1]{MR}.

\begin{cor}\label{cor:NoPerc}
Let $\SP=\R^d$ and $\l$ be Lebesgue measure.
If $$z(J):=\limsup_{|x|\to \infty} J(x,x)<z_c$$ then
$\m(\,\exists \text{ infinite cluster of }B_R(\cdot)\,)=0$.
\end{cor}

We now address the question of whether DPP's are Gibbsian. In view of
Remark \ref{rem:MWM} and Theorem \ref{thm:MWM}, the following
result implies that $\m$ is Gibbsian at least in a general sense.

\begin{thm}\label{thm:Sigma}
$\m$ satisfies condition $(\S_\l)$.  Its CPI is given by
\begin{equation}\label{eq:cJ-lim}
\hat c(\a,\x) = \lim_{n\to\infty} \hat c_{\D_n}(\a,\x_{\D_n})\quad 
\text{for $L\otimes\m$-almost all $(\a,\x)$,}
\end{equation}
where $\hat c_{\D_n}$ is given by \eqref{eq:cLformula} 
and $(\D_n)$ is any
sequence in $\cB_0$ that increases to $\SP$. In particular, 
$\m(N_\L{=}0\,|\,\cF_{\L^c}) >0$ $\m$-almost surely for
each $\L\in\cB_0$ and, for any $\cF_\L$-measurable $f:\cX\to\R_+$, the ratio
\[
\E_{\m}(f\,|\,\cF_{\L^c})/\m(N_\L{=}0\,|\,\cF_{\L^c})
\]
of conditional expectations has a decreasing version.
\end{thm}
Two remarks are in order. First, in view of the comments in 
Remark \ref{rem:G_L,D}(b) and below \eqref{eq:Papint},
the first assertion can be restated as an absolute continuity 
between DPP's: the DPP $\m^x$ with correlation
kernel $K^x$ is absolutely continuous with respect to the 
DPP $\m$ for $K$ with density $c(x, \cdot)/K(x,x)$.
Next, and more importantly, combining  \eqref{eq:cJ-lim} 
with \eqref{eq:cLformula} and \eqref{eq:condProb} we
obtain at least an implicit formula for the conditional probabilities 
of $\m$ given the events outside of a
local region. But the question remains of whether the CPI, and 
thereby the conditional probabilities, of $\m$
can be identified in a more specific way. In fact, if we make 
a proper choice of the integral kernel $J(x,y)$ as
in Lemma \ref{lem:loc_trace_formula}, a natural candidate for $\hat c$ is
\begin{equation}\label{eq:candidate}
\hat c_*(\a,\x) = \lim_{\D\uparrow\SP} \det J(\a\x_\D,\a\x_\D)/\det J(\x_\D,\x_\D)\,.
\end{equation}
This limit does exist because the ratio on
the right-hand side is decreasing in $\D$, as follows from inequality
\eqref{eq:det3} in the Appendix. (As before, we set a ratio of determinants
equal to zero if the denominator vanishes.) In contrast to
\eqref{eq:cJ-lim}, the determinants in \eqref{eq:candidate}
involve $J$ itself rather than $J_{[\D_n]}$.
Now our question is the following: \emph{When is it true that
$\hat c=\hat c_*$ a.e. for $L\otimes\m$?}

Unfortunately, we are unable to settle the question above 
in the same generality as this was done in the lattice
case by Shirai and Takahashi \cite[Theorem 6.2]{ST2}. 
Their argument exploits the symmetry between occupied and
empty lattice sites and therefore does not carry over to our 
continuous setting. (Formally, this is reflected by
their condition $\text{spec\,} K\sub\; ]0,1[$ which is never 
satisfied in the continuous case.) The following
theorem gives a partial answer, at least.

\begin{thm} \label{thm:Gibbs}
In general, the inequality $\hat c\le \hat c_*$ holds
$L\otimes\m$-almost everywhere. If $\hat c_*$ is continuous
$L\otimes\m$-a.e. then $\hat c= \hat c_*$
$L\otimes\m$-a.e., i.e., $\hat c_*$ is a version
of the CPI of $\m$.
\end{thm}

Proposition \ref{prop:a.e.continuity} below will provide
sufficient conditions for the almost-everwhere continuity
to hold. First let us comment  on the significance of this result.

\begin{rem}\label{rem:potential}
If $\hat c= \hat c_*$ a.e. then $\m$ is a
Gibbs measure for a specification $G$ defined in terms of $J$
as follows. For each $\L\in\cB_0$ and any
$\x\in\cX_{\L^c}$ satisfying $\det J(\z,\z)>0$ for all finite
$\z\sub\x$, let $G_\L(\cdot\,|\,\x)$ be defined by
inserting $\hat c _*$ into \eqref{eq:condProb}.  In fact,
$G_\L(\cdot\,|\,\x)$ is a Gibbs distribution for the Hamiltonian
$H_\L(\cdot\,|\,\x)= -\log \hat c _*(\cdot,\x)$ on $\cX_\L$,
$\x\in\cX_{\L^c}$. (If desired, one can express the Hamiltonian in
terms of a many-body potential $\Phi$, but this is not
particularly useful.) These Gibbs distributions altogether form a
Gibbsian specification $G=(G_\L)_{\L\in\cB_0}$ in the sense of
Preston \cite[pp. 16, 17]{PrLNM}; this has been proved by Gl\"otzl
\cite{Gloetzl2} in a general setting, and in \cite{Y} for the
particular case of DPP's. In particular, if $\hat c_*$ is a.e.
continuous then the conditional expectations
$\E_\m(g|\cF_{\L^c})$ in \eqref{eq:condexpect} have a.e.\ continuous
versions. This is Gibbsianness in the best sense one
can expect for a point process.

Conversely, suppose $\m $ is a Gibbs measure for some Hamiltonian
$H$. Then $\hat c^H(\a,\x)$ $:=\exp[-H_\a(\a\,|\,\x_{\a^c})]$ is a
version of its CPI \cite{NZ} and will typically satisfy the continuity
condition $\hat c^H(\a,\x) = \lim_{\D\ua\SP}\hat c^H(\a,\x_\D)$.
If one assumes that the function $\hat c =\lim_{n\to\infty} \hat c
_{\D_n}$ in  \eqref{eq:cJ-lim} has the same continuity property,
one can conclude that $\hat c =\hat c _*$ a.e. For, Proposition
\ref{prop:local_convergence} below implies that $\lim_{n\to\infty} \hat
c _{\D_n}(\a,\x_\D) =\hat c _*(\a,\x_\D)$ a.e.
for all $\D\in\cB_0$.
\end{rem}

The next proposition describes a case in which the condition of
almost-everwhere continuity in Theorem \ref{thm:Gibbs} is satisfied.

\begin{prop}\label{prop:a.e.continuity}
Suppose $E=\R^d$, $\l$ is Lebesgue measure and, in addition to (H),
\begin{enumerate}
\item $J:=K(I-K)^{-1}$ has a continuous integral kernel $J(x,y)$, $x,y\in\R^d$;
\item $J$ has finite range $R<\infty$, in that $J(x,y)=0$ for $|x-y|>R$;
\item $\m$ does not percolate, in that
$\m \big(\,\exists \text{ \rm infinite cluster of }B_R(\cdot) \,\big)=0$.
\end{enumerate}
Then $\hat c_*$ is continuous $L\otimes\m$-almost everywhere.
Specifically,
\begin{equation}\label{eq:cJ-subcr}
\hat c _*(\a,\x)=\det J(\a\x_W,\a\x_W)/ \det
J(\x_W,\x_W)\ 1_{\{\diam W(\a,\x)<\infty\}}\,,
\end{equation}
where $W(\a,\x)$ is the
union of the clusters of $B_R(\a\x)$ hitting $\a$, and
$\x_W:=\x_{W(\a,\x)}$ is the restriction of $\x$ to $W(\a,\x)$.
\end{prop}

Under the conditions of this proposition, the Gibbsianness of DPP's has
already been derived by the second author in \cite{Y} by different
methods.

\begin{ex}\label{ex:Example}
Typical examples of operators satisfying the assumptions
of Proposition \ref{prop:a.e.continuity} are
obtained by Fourier transforms: Let $\hat j\ge0$ be an integrable
function on $\R^d$,
\[
j(x):=(2\p)^{-d} \int_{\R^d}e^{ix\cdot t}\;\hat j(t)\,dt
\]
its (inverse) Fourier transform, and $J(x,y):=j(x-y)$. We assume
that $j(\cdot)\in L^1(\R^d)$. By Bochner's theorem and Young's
inequality \cite[p. 29]{RS}, $J$ then defines a positive bounded
operator on $L^2(\R^d)$. The associated  correlation operator
$K:=J(I+J)^{-1}$ is the convolution operator for the function $k$
satisfying $\hat k = \hat j/(1+\hat j)$. The validity of
hypothesis (H) is evident,
and the continuity assumption (a) of Proposition
\ref{prop:a.e.continuity} also holds by definition.  To satisfy
the finite range condition (b), let $j:=\ph\,\chi_R$ or,
equivalently, $\hat j=\hat \ph\star \hat \chi_R$, where $0\le
\hat\ph\in L^1(\R^d)$, $\chi_R(x)=\prod_{i=1}^d(1-|x_i|/R)^+$, and
$\star$ denotes convolution. By Corollary \ref{cor:NoPerc}, the
non-percolation assumption (c) holds  when $d=1$ or
$j(0)<z_c(d,R)$.

To construct non-translation-invariant examples, let $\psi\ge0$ be any
bounded measurable function on $\R^d$
and $M_\psi$ the associated multiplication operator on $L^2(\R^d,\l)$.
Let $J$ be a finite range operator as above
and $L:=JM_\psi J$. $L$ has the kernel $L(x,y)=\int
 j(x-z)\psi(z)j(z-y)\,dz$. Obviously, $L(x,y)$ is of finite range
because so is $j$, and not translation-invariant except when $\psi$ is
constant almost everywhere.
\end{ex}

Let us also comment on the particular case of renewal processes 
which deserves special interest.

\begin{ex}\label{ex:renewal}
Let $\SP=\R$ and $K:f\to k\star f$ the convolution operator for the function
$k(x)= \r e^{-a|x|}$, where $\r,a>0$ are such that $\r <a/2$. (The last
condition means that $\|k\|_1<1$, and thus $\|K\|<1$ by Young's inequality.)
In the setting of the preceding example, this corresponds to setting
$\hat k(t) = 2\r a/(a^2+t^2)$. Hence
\[
\hat j(t) = \frac{\hat k(t)}{1-\hat k(t)}
=\frac{2\r a}{\s^2+t^2}
\]
with $\s^2:= a^2-2\r a$, and therefore $j(x)=(\r a/\s) e^{-\s|x|}$. The
associated integral kernel can be written in the form
\[
J(x,y):=j(x-y) = u(x\wedge y)\, v(x\vee y)
\]
with $x\wedge y=\min(x,y)$, $x\vee y=\max(x,y)$, $u(x)= e^{\s x}$, and
$v(x)= (\r a/\s) e^{-\s x}$. Therefore, if $\a=\{x_1,\ldots,x_n\}$ with 
$x_1<\cdots<x_n$ then
\[
\det J(\a,\a) = u(x_1)\, v(x_n)\, \prod_{i=1}^{n-1}d(x_{i+1}-x_i)
\]
with $d(x_{i+1}-x_i)= \det \begin{pmatrix}v(x_i)&v(x_{i+1})\\
u(x_i)&u(x_{i+1}) \end{pmatrix} = (2\r a/\s) \sinh\big(\s(x_{i+1}-x_i)\big)$, as
is easily verified by induction on $n=|\a|$. Together with \eqref{eq:candidate}
we find that
\[
c_*(x,\x) = {d\big(\ell_x(\x)\big)\, d\big(r_x(\x)\big)}\big/
{d\big(\ell_x(\x)+r_x(\x)\big)}\,,
\]
where $\ell_x(\x)$ and $r_x(\x)$ are the distances from $x$ to the nearest
points of $\x$ on the left, respectively on the right. This relationship still
holds if $d(\cdot)$ is multiplied with an exponential factor to obtain a probability
density. That is, we have
\[
c_*(x,\x) = {f\big(\ell_x(\x)\big)\, f\big(r_x(\x)\big)}\big/{f\big(\ell_x(\x)+r_x(\x)\big)}
\]
for the probability density $f(x) = e^{-ax}d(x)$ on  $[0,\infty[$.
The right-hand side is the PI of the stationary renewal process $\m_f$
with spacing density $f$, and $\m_f$ is the only translation invariant 
point process having this PI. As $c_*$ is continuous,
we conclude from Theorem \ref{thm:Gibbs} that
$\m_f$ is the unique DPP for $K$. In this way we recover a well-known 
result of Macchi \cite{M}; see also
\cite{DV}. According to Soshnikov \cite{So}, the densities $f$ of the 
form above are the only possible spacing
distributions of determinantal renewal processes. (We remark that 
Soshnikov also gave a complete
characterization of all kernels $K(x,y)$ on $\SP=[0,\infty[$ which 
are such that the associated DPP's have
independent spacings, see also Example~6.7 of \cite{ST1}. 
The previous comments can be extended to this setting.)
\end{ex}

\section{Proofs}
\label{sec:proofs}

Let us start with the proof of Theorem \ref{thm:PapMon}. Its key
ingredients are equation \eqref{eq:density_function} for the
Janossy densities of $\m$, and the determinant inequalities of
Lemma \ref{lem:determinant-relation} in the Appendix.
We assume throughout this section that hypothesis (H) is satisfied.

\medskip
\Proof[ of Theorem \ref{thm:PapMon}]  Fix any $\L\in\cB_0$.
>From Eq. (\ref{eq:density_function}) we know that
$\m_\L$ is absolutely continuous relative to $L_\L$ with
density $\s_\L(\x)\propto \det J_{[\L]}(\x,\x)$.
Inequality \eqref{eq:det2} implies that the zero
set $\{\s_\L=0\}$ is increasing. Remark \ref{rem:locMWM} thus
tells us that $\m_\L$ satisfies condition $(\S_{\l_\L})$ with
PI given by \eqref{eq:cLformula}. This proves the
first statement of the theorem.

For the proof of \eqref{eq:cLinequality} let $\x\sub\eta\in \cX_\L$
and $\a\in\cX_0$ with $\a\sub \L\sm\eta$.
Choosing the kernel function of $J_{[\L]}$ as in Lemma
\ref{lem:positive_kernel}, we get that $J_{[\L]}(\a\eta,\a\eta)$ is
positive definite. Inequality \eqref{eq:det3} thus asserts that
\[
\det J_{[\L]}(\a\eta,\a\eta)\det J_{[\L]}(\x,\x) \le \det
J_{[\L]}(\eta,\eta) \det J_{[\L]}(\a\x,\a\x),
\]
and this gives the
first inequality in \eqref{eq:cLinequality}. In particular, we
have $\hat c_\L(\a,\x)\le c_\L(\a,\emptyset)$ for the empty configuration
$\emptyset$. But the last quantity equals $\det
J_{[\L]}(\a,\a)$ by definition. An application of inequality \eqref{eq:det2}
then gives the last estimate.\EndProof

To exploit the upper bound in \eqref{eq:cLinequality} we will
use the following comparison lemma.

\begin{lem}\label{lem:monotonicity}
For any $\L \subset \D\in \cB_0$,
$J_{[\L]}\le P_\L J_{[\D]}P_\L\le  J_\L$ in the operator ordering.
In particular, 
\begin{equation}\label{eq:det-monotone}
\det J_{[\L]}(\x,\x)  \le 
\det J_{[\D]}(\x,\x)\le \det J_\L(\x,\x)
\end{equation}
for $L_\L$-almost all $\x\in \cX_\L$.
\end{lem}

\Proof We prove only the inequality $J_{[\L]}\le J_\L$ which
implies in particular that $P_\L J_{[\D]}P_\L\le P_\L J_\D P_\L
=J_\L$. The proof of the relation $J_{[\L]}\le P_\L J_{[\D]}P_\L$
is similar. Since $J=K (I-K )^{-1}$ and
$L(I-L)^{-1}=-I+(I-L)^{-1}$ for any contraction $L$, the
inequality $J_{[\L]}\le J_\L$ is equivalent to
\[
P_\L(I-K_\L )^{-1}P_\L\le P_\L(I-K )^{-1}P_\L.
\]
But this follows from Lemma \ref{lem:basic-order} as applied to
$P=P_\L$ and $T=I-K $, because $P_\L(I-K_\L
)^{-1}P_\L=P_\L(P_\L(I-K )P_\L)^{-1}P_\L$.

To prove the determinantal inequalities \eqref{eq:det-monotone},
we apply Lemma \ref{lem:positive_kernel} to the positive trace
class operators $J_{[\L]}$, $P_\L J_\D P_\L-J_{[\L]}$, and
$J_\L-P_\L J_\D P_\L$. We then see that the kernels of $J_{[\L]}$,
$J_\D$, and $J_\L$ can be chosen in such a way that, for each
$\x\in \cX_\L$, $J_{[\L]}(\x,\x)\le  J_\D(\x,\x)\le J_\L(\x,\x)$
as operators on $\C^{|\x|}$. By Lemma \ref{lem:indistinguishable},
these kernels are indistinguishable from the corresponding kernels
obtained by applying Lemma \ref{lem:positive_kernel} to
$J_{[\L]}$, $J_\D$, and $J_\L$ directly. For the latter kernels,
the same operator inequalities thus hold for $L_\L$-almost all
$\x\in \cX_\L$, at least. To complete the proof it is therefore
sufficient to note that the determinant is increasing relative to
the operator ordering; see \cite[Corollary III.2.3]{B}.\EndProof

We are now ready for the proofs of Remark~\ref{rem:G_L,D} and
Corollaries~\ref{cor:Monoton} to \ref{cor:NoPerc}.

\medskip
\Proof[ of Remark \ref{rem:G_L,D}] (a) Let us start with a
slightly simpler estimate. Combining  the last bound in
\eqref{eq:cLinequality} with the second inequality in
\eqref{eq:det-monotone} (for $\x=\{x\}$) and using \eqref{eq:L}
and the trace formula \eqref{eq:trace-formula} we find
\[
Z_{\L,\D}(\x)\le \exp\int_\L J_\L(x,x)\,\l(dx)= \exp \Tr J_\L <\infty
\]
because $J_\L$ is trace class by (H).

To obtain the sharper bound stated in Remark \ref{rem:G_L,D}(a)
we use the next to last inequality in \eqref{eq:cLinequality} and
the second inequality in \eqref{eq:det-monotone} (for $\x=\a$)
to obtain
\[
Z_{\L,\D}(\x)\le \int L_\L(d\a)\,\det J_\L(\a,\a)= \det(I+J_\L).
\]
To get the last identity one can use the argument on p. 930 of \cite{So}
by combining Lidskii's theorem \cite[p. 50]{Si} with the trace formula
\eqref{eq:trace-formula}. (A similar argument implies that the Janossy densities
in \eqref{eq:density_function} integrate to $1$.)
We note that $\det(I+J_\L)\le \exp\Tr (J_\L)$ by \cite[eq. (3.6)]{Si}.

(b) For fixed $\x$ as stated, $G_{\L,\D}(\,\cdot\,|\,\x)$ -- considered as a measure
on $\cX_\L$ -- is defined by an $L_\L$-density of the form
\[
\s^\x_{\L,\D}(\a)\propto \hat c_\D(\a, \x) = \det
J_{[\D]}(\a\x,\a\x)/\det J_{[\D]}(\x,\x)\,,\quad \a\in\cX_\L.
\]
In view of the identity \eqref{eq:det-ratio}, the last expression is equal to
$\det J_{\L,\D}^\x(\a,\a)$. Comparing this observation with
\eqref{eq:density_function} we see that $G_{\L,\D}(\,\cdot\,|\,\x)$
is a determinantal process with interaction kernel $J_{\L,\D}^\x$,
defining a positive Hermitian operator on $L^2(\L,\l_\L)$. The associated
correlation operator is
$K_{\L,\D}^\x:=J_{\L,\D}^\x(I+J_{\L,\D}^\x)^{-1}$.\EndProof

\medskip
\Proof[ of Corollary \ref{cor:Monoton}] In view of Eq. (\ref{eq:locCondProb}),
\[
G_{\L,\D}(f\,|\,\x )\,\big/\, G_{\L,\D}(N_\L{=}0\,|\,\x )= \int L_\L(d\a)\,
f(\a) \, \hat c_\D(\a,\x)
\]
for $\m_\D$-almost all $\x\in\cX_{\D\sm\L}$. Since $f\ge0$,
Theorem \ref{thm:PapMon}(b) together with equation
\eqref{eq:compPapdens} implies that the integrand on the
right-hand side is a decreasing function of $\x$. In particular,
for $f\equiv1$ we find that $G_{\L,\D}(N_\L{=}0\,|\,\cdot)$ is
increasing. To prove \eqref{eq:Monoton} we note that
\[
\begin{split}
&\m(N_\L{=}0\,|\,N_{\D_2\sm\L}{=}0)= G_{\L,\D_2}(N_\L{=}0\,|\,\emptyset )\\
&\le
\E_{\m}\big(G_{\L,\D_2}(N_\L{=}0\,|\,\cdot)\big|\cF_{\D_1\sm\L}\big)(\emptyset)\\
&= \m(N_\L{=}0\,|\,\cF_{\D_1\sm\L})(\emptyset)
=\m(N_\L{=}0\,|\,N_{\D_1\sm\L}{=}0)
\end{split}
\]
because $G_{\L,\D_2}(N_\L{=}0\,|\,\emptyset )\le
G_{\L,\D_2}(N_\L{=}0\,|\,\cdot )$. The final inequality follows by
setting $\D_1=\L$.\EndProof

\medskip
\Proof[ of Corollary \ref{cor:PoissonDom}] We show first that
$\m_\L\preceq\p^J_\L$ for any $\L\in\cB_0$. Consider the PI $c_\L$
of $\m_\L$. Combining \eqref{eq:cLinequality} (for $\a=\{x\}$)
with \eqref{eq:det-monotone} (for $\x=\{x\}$) we find that
$c_\L(x,\x)\le J_\L(x,x)$ for all $x\in \L$ and $\x\in \cX_\L$. 
On the other hand, it is
well-known and easy to check that the Poisson point process
$\p^J_\L$ has the PI $c_\L^J(x,\x)=J(x,x)$ on $\L$; cf.
\cite{Mecke, NZ}. Moreover, since $J(x,x)=J_\L(x,x)$ for $\l$-a.a.
$x\in \L$ by the construction of $J(x,x)$ in Lemma
\ref{lem:loc_trace_formula}, we may also regard $J_\L(x,x)$ as the
PI of $\p^J_\L$. Thus we can conclude that the PI of $\m_\L$ at a
configuration $\x\in \cX_\L$ is not larger than the PI of
$\p^J_\L$ evaluated at any larger configuration $\eta\in\cX_\L$,
as long as $x\notin\eta$. This, however, is precisely the
hypothesis of the point-process counterpart of the well-known
FKG-Holley-Preston inequality for lattice systems. This continuous
counterpart was obtained first by Preston \cite{Pr}; an
alternative simplified proof can be found in \cite[Theorem
1.1]{GK}. It asserts that, under the above condition on the PI's,
$\m_\L\preceq\p^J_\L$, as required.

To get rid of the locality restriction we argue as follows. For
any compact $\L$, a celebrated theorem of Strassen \cite{Strassen}
provides us with a probability measure $m_\L$ on $\cX\times\cX$
having marginals $\m_\L$ resp. $\p^J_\L$ and being supported on
the set $D:=\{(\x,\eta)\in\cX\times\cX:\x\sub\eta\}$. Note that
$D$ is closed when $\cX\times\cX$ is equipped with the product of
the vague topology on $\cX$. By a standard compactness criterion
for point processes \cite[Proposition 9.1.V]{DV}, the measures
$m_\L$ admit a weak limiting measure $m$ as $\L$ increases to
$\SP$. By construction, $m$ has marginals $\m$ resp. $\p^J$, and
$m(D)=1$ because $D$ is closed. This implies that
$\m\preceq\p^J$.\EndProof

\Proof[ of Corollary \ref{cor:NoPerc}] By assumption, there exists
some $z<z_c$ such that $J(x,x)\le z$ for all $x$ outside of some
compact region. Therefore, $\p^J$-almost surely there exists no
infinite cluster outside of this region, and therefore no infinite
cluster anywhere. As the existence of an infinite cluster is an
increasing event, the absence of percolation carries over to $\m$,
by Corollary \ref{cor:PoissonDom}.\EndProof

Next we turn to the proof of Theorem~\ref{thm:Sigma}. We will use
a martingale argument to derive property $(\S_\l)$ from the local
property $(\S_{\l_\L})$ of Theorem~\ref{thm:PapMon}.

\medskip
\Proof[ of Theorem~\ref{thm:Sigma}] We only need to show
that $\m$ satisfies condition $(\S_\l)$ with CPI
\eqref{eq:cJ-lim}; the remaining assertions then follow
from Theorem \ref{thm:PapMon} and Remark
\ref{rem:MWM}(c). Let $(\D_n)$ be any increasing sequence
in $\cB_0$ exhausting $\SP$, $\L\in\cB_0$ a fixed set,
and $n$ so large that $\L\sub\D_n$. Consider the product
space $\cX_\L\times\cX$, equipped with the probability
measure $\n_\L:=\p_\L^1\otimes\m$ and the $\s$-algebras
$\cG_n= (\cF|_{\cX_\L})\otimes\cF_{\D_n}$. Now, Theorem
\ref{thm:PapMon} asserts that, on $\cG_n$, the restriction of
$\CCM[\m]$ to $\cX_\L\times\cX$ is absolutely
continuous with respect to $\n_\L$ with Radon-Nikodym
density $R_n:=e^{\l(\L)}\,\hat c_{\D_n}$.
The sequence
$(R_n)$ is therefore a nonnegative martingale relative to
$\n_\L$. By \eqref{eq:cLinequality} and
\eqref{eq:det-monotone},
$R_n$ satisfies the bound
\[
R_n(\a,\x)\le e^{\l(\L)}\,\det J_{[\D_n]}(\a,\a)
\le S(\a):= e^{\l(\L)}\,\det J_{\L}(\a,\a)
\]
for $\n_\L$-almost all $(\a,\x)\in\cX_\L\times\cX$.
As we have seen in the proof of Remark \ref{rem:G_L,D}(a),
$S$ is integrable relative to $\p_\L^1$.
As a consequence, $(R_n)$ converges $\n_\L$-almost surely and 
in $L^1(\n_\L)$-norm to a limit $R$. By the
norm-convergence, $R$ is a Radon-Nikodym density of $\CCM[\m]$ 
relative to $\n_\L$ on the limiting $\s$-algebra
$\s(\bigcup_n\cG_n)= (\cF|_{\cX_\L})\otimes\cF$. Finally we 
replace $\p_\L^1$ with $L_\L$ by dropping the
constant factor $e^{\l(\L)}$,  and use the fact that 
$\cX_\L\uparrow\cX_0$ as $\L \uparrow\SP$. We then get the
desired result that $\CCM[\m]\ll L\otimes\m$ with density 
\eqref{eq:cJ-lim}.\EndProof

Our proof of Theorem \ref{thm:Gibbs} will be based on the
following convergence result. As before, we say that a kernel
function $t(x,y)$ of an operator $T$ is \emph{positive definite
almost everywhere} if the matrix $t(\a,\a)$ is positive definite
for $L$-almost all $\a\in \cX_0$.

\begin{prop}\label{prop:local_convergence}
Given any sequence  $(\D_n)$ in $\cB_0$ increasing to $\SP$, there
exist integral kernels $J_{[n]}(x,y)$ for $J_{[n]}:=J_{[\D_n]}$
such that the following holds:
\begin{enumerate}
\item Each $J_{[n]}(x,y)$ is positive definite a.e. and satisfies
the trace formula
\begin{equation}\label{eq:locTrace2}
\Tr\,(P_\L J_{[n]}P_\L)=\int_\L J_{[n]}(x,x)\,\l(dx)
\end{equation}
for all $\L\subset\D_n$. \item The limit
\[
J(x,y)=\lim_{n\to \infty} J_{[n]}(x,y)
\]
exists for all $x,y\in\SP\sm\cN$, where $\cN\in \cB$ is a
$\l$-nullset. Moreover, this limit defines an integral kernel of
$J$, is positive definite a.e., and satisfies the trace formula
\eqref{eq:locTrace1}.
\end{enumerate}
\end{prop}
Note that, by Lemma \ref{lem:indistinguishable}, the kernels
$J_{[n]}(x,y)$ and $J(x,y)$ constructed here are indistinguishable
from those obtained directly from Lemmas \ref{lem:positive_kernel}
and \ref{lem:loc_trace_formula}.

\bigskip
\Proof
In view of hypothesis (H) and Lemma \ref{lem:loc_trace_formula},
$K$ admits an integral kernel $K(x,y)$ with the properties
(a)--(c) of that lemma. In particular, $K(x,y)$ is positive
definite a.e., and there exists a $\l$-nullset $\cN\in \cB$ such
that $k_x:=K(\cdot,x)\in L^2(\SP,\l)$ whenever $x\notin \cN$. For
each $n$ we choose the integral kernel
$K_n(x,y):=1_{\D_n}(x)K(x,y)1_{\D_n}(y)$, $x,y\in\D_n$, of
$K_n:=K_{\D_n}$.  Also, for each power $l\ge 2$ we set
$K^l(x,y):=(k_x,K^{l-2} k_y)$ and $K_n^l(x,y):=(k_x,P_n(K_n)^{l-2}
k_y)$, where $(\cdot,\cdot)$ is the inner product in $L^2(\SP,\l)$
and $P_n:=P_{\D_n}$. These are well-defined for $x,y\notin \cN$,
and define versions of the integral kernels of $K^l$ resp.
$(K_n)^l$.

Now, since $\|K_{n}\|\le \|K\|<1$, we can write
$J_{[n]}=\sum_{l\ge1}(K_n)^l$ and $J=\sum_{l\ge1}K^l$\,; these
series converge in operator norm, and the convergence
is uniform in $n$. In particular, this implies the
absolute convergence of the series $\sum_{l\ge 1}K_n^l(x,y)$
and $\sum_{l\ge 1}K^l(x,y)$ for $x,y\notin \cN$,
again uniformly in $n$. Obviously, these series define
integral kernels of $J_{[n]}$ and $J$, respectively. They
are positive definite a.e. because $K(x,y)$ is.

To prove the trace formula \eqref{eq:locTrace1} we introduce the
partial sum $J^{(m)}=\sum_{l=1}^mK^l $.
Obviously, as $m\to \infty$, the sequence $P_\L J^{(m)}P_\L$
increases and converges in operator norm to $J_\L$, which is of
trace class. Using \cite[Theorem 2.16]{Si}, we conclude that $P_\L
J^{(m)}P_\L\to J_\L$ in trace norm as $m\to \infty$. Hence
\[
\Tr\big(J_\L\big)=\lim_{m\to \infty}\,\sum_{l=1}^m\Tr\big(P_\L K^l
P_\L\big)\,.
\]
Inserting the trace formulas of Lemma \ref{lem:loc_trace_formula}(c)
and taking the limit $m\to \infty$ we arrive at
\eqref{eq:locTrace1}. Formula \eqref{eq:locTrace2} follows in the same way by considering $J_{[n]}$ resp.
$K_n$ instead of $J$ resp. $K$.

It remains to show that $J_{[n]}(x,y)\to J(x,y)$ when $x,y\notin\cN$.
By the uniform absolute convergence of these series, it is enough
to check that $K_{n}^l(x,y)\to K^l(x,y)$ for all $l\in\N$.
The case $l=1$ is trivial. For  $l=2$ we have
$K_{ n}^2(x,y)=(k_x,P_{n} k_y)$
when $\D_n\ni x,y$, and this converges to $(k_x,
k_y)=K^2(x,y)$. In the case $l\ge 3$ we note that $K_n$
converges strongly to $K$ as $n\to \infty$. Since $\|K_n\|\le1$,
it follows that $(K_{n})^{l-2}$ converges strongly to $K^{l-2}$,
whence $K_{n}^l(x,y)=(k_x,(K_{n})^{l-2}k_y)$
converges to $(k_x,K^{l-2}k_y)=K^l(x,y)$. This
completes the proof.\EndProof

\Proof[ of Theorem \ref{thm:Gibbs}] Consider the CPI
\[
\hat c (\a,\x) = \lim_{n\to\infty} \hat c _{\D_n}(\a,\x_{\D_n})\,.
\]
that occurs in \eqref{eq:cJ-lim}. We know that this limit exists
for $L\otimes\m $-almost all $(\a,\x)\in\cX_0\times\cX$. In fact,
we may define $\hat c _{\D_n}(\a,\x_{\D_n})$ in terms of the
integral kernels $J_{[n]}$ of Proposition
\ref{prop:local_convergence} by setting $\hat c
_{\D_n}(\a,\x_{\D_n})=\det J_{[n]}(\a\x_{\D_n},\a\x_{\D_n}) / \det
J_{[n]}(\x_{\D_n},\x_{\D_n})$. On the other hand, for each
$\L\in\cB_0$ we conclude from \eqref{eq:cLinequality} that $\hat
c_{\D_n}(\a,\x_{\D_n}) \le \hat c _{\D_n}(\a,\x_{\L})$ for
$L\otimes \m$-a.a. $(\a,\x)$ as soon as $\D_n\supset\L$ . Hence
\[
\hat c (\a,\x) \le \liminf_{n\to\infty} \hat c_{\D_n}(\a,\x_\L)\,.
\]
Now, Proposition \ref{prop:local_convergence} implies that for
$L\otimes L_\L$-a.a. $(\a,\x_\L)$
\[
\hat c _{\D_n}(\a,\x_\L)\to \hat c _*(\a,\x_\L)= \det
J(\a\x_\L,\a\x_\L)/\det J(\x_\L,\x_\L)
\]
as $n\to\infty$, provided that $\det J(\x_\L,\x_\L)>0$. In view of
\eqref{eq:density_function} and \eqref{eq:det-monotone}, the last
proviso holds true for $\m$-almost all $\x$. Hence $\hat c
(\a,\x) \le \hat c _*(\a,\x_\L)$ for $L\otimes\m$-almost all
$(\a,\x)$. As we have noticed after \eqref{eq:candidate},
$\hat c_*(\a,\x_\L)$ is a decreasing function of $\L$, and its limit is
equal to $\hat c _*(\a,\x)$ by definition. This proves that
$\hat c \le \hat c_*$ almost everywhere.

Suppose next that $\hat c_*$ is continuous
$L\otimes\m$-everywhere. Consider the DPP's $\m^{J_\D}$ for the
interaction operators $J_\D=P_\D J P_\D$, and thus for the
correlation operators $K_{[\D]}=J_\D(I+J_\D)^{-1}$,
$\L\sub\D\in\cB_0$. According to Theorem 6.17 of \cite{ST1}, or
Lemma 4.3 of \cite{Y} in the case of a discontinuous kernel of
$J$, $\m^{J_\D}$ converges weakly to $\m $ as $\D$ increases to
$\SP$. Let $f$ be any bounded continuous function on
$\cX_0\times\cX$ such that $f(\a,\x)=0$ unless $\a\sub\L$ and
$|\a\x_\L|\le k$ for some number $k$. The portmanteau theorem
\cite[Proposition A2.3.V]{DV} then implies that 
\[
\int f \,\hat  c _*\, d(L_\L\otimes\m ) = \lim_{\D\ua\SP} \int f
\,\hat  c _* \,d (L_\L\otimes\m^{J_\D})\,.
\]
But $\m^{J_\D}$ is supported on $\cX_\D$, and for almost every
$\a\in \cX_\L$ and $\x\in \cX_\D$ we have  $\hat c_*(\a,\x)=\hat
c^{J_\D}_*(\a,\x):=\det J_\D(\a\x,\a\x)/\det J_\D(\x,\x)$. By
Remark \ref{rem:locMWM},  $\hat c^{J_\D}_*$ is a CPI of
$\m^{J_\D}$. Eq. \eqref{eq:compPapint} thus shows that
\[
\int f \,\hat  c _*\,d (L_\L\otimes\m^{J_\D}) = \int \m^{J_\D}(d\x
) \sum_{\a\in\cX_0:\,\a\subset \x } f(\a, \x \sm \a)\,.
\]
The integrand on the right-hand side is still a bounded continuous
function of $\x$. Letting $\D\ua\SP$ we thus find that
\[
\int  f \,\hat  c _*\, d(L_\L\otimes\m )= \int  f\, d\CCM[\m
]|_{\cX_\L\times \cX}\,.
\]
As $f$, $k$ and $\L$ were arbitrarily chosen, it follows that
$\hat c _*$ is a CPI of $\m $.\EndProof

Finally we prove Proposition \ref{prop:a.e.continuity} on the a.e.
continuity of $\hat c_*$.

\medskip
\Proof[ of Proposition \ref{prop:a.e.continuity}] Fix any $\L\in\cB_0$,
$\a\in\cX_\L$, and consider the function $\hat c_*(\a,\cdot)$. For
$\x\in \cX$ let $W(\a,\x)$ be the union of all clusters of
$B_R(\a\x)$ hitting $\a$, and write $\x_W=\x_{W(\a,\x)}$. By the
finite range assumption,  $J(x,y)=0$ whenever $x\in \a\x_W$ and
$y\in \x\sm \x_W$. Hence, if $W(\a,\x)$ is finite then, for any
local $\D\supset W(\a,\x)$, the matrix $J(\a\x_\D,\a\x_\D)$
consists of the two diagonal blocks $J(\a\x_W,\a\x_W)$ and
$J(\x_{\D\sm W},\x_{\D\sm W})$. Thus
\[
\det J(\a\x_\D,\a\x_\D)= \det J(\a\x_W,\a\x_W)\, \det J(\x_{\D\sm
W},\x_{\D\sm W})\,.
\]
Dividing this by the analogous equation for $\a=\emptyset$ we
arrive at \eqref{eq:cJ-subcr}. By the continuity of the kernel $J$,
it follows that $\hat c_*(\a,\cdot)$ is continuous on the open
set of configurations having no infinite cluster, which has full measure
by assumption.\EndProof

\renewcommand{\thesection}{\Alph{section}}
\setcounter{section}{0}
\section{Appendix}
Here we collect some general facts on positive matrices and
operators. Some of them are standard and listed here for ease of
reference. (We include their simple proofs for the convenience of
the readers.) Some others are more subtle and do not seem to
appear in the literature, at least not in the form we need. We
start with a classical lemma on determinants of finite matrices.
The inequalities \eqref{eq:det2} and \eqref{eq:det3} are a key ingredient
of our arguments; \eqref{eq:det2} is known as Fischer's inequality. To state
the lemma we consider a finite index set $\G$ and any complex matrix
$A=(a_{ij})_{i,j\in \G}$. For any two subsets $\a,\b\sub \G$ we let
$A_{\a,\b}=(a_{ij})_{i\in \a,j\in \b}$ be the associated submatrix.

\begin{lem}\label{lem:determinant-relation}
If $A$ is Hermitian and positive definite then
\begin{equation}\label{eq:det2}
\det A \le \det A_{\a,\a} \,\det A_{\b,\b}
\end{equation}
when $\G=\a\cup\b$ is partitioned into two nontrivial subsets, and
\begin{equation}\label{eq:det3}
\det A \;\det  A_{\b,\b} \le \det A_{\a\cup\b,\a\cup\b}
 \,\det A_{\b\cup\g,\b\cup\g}
\end{equation}
when $\G=\a\cup\b\cup\g$ is partitioned into three nontrivial subsets.
For $\G=\a\cup\b$  and general $A$ with invertible $A_{\b,\b}$ we have
\begin{equation}\label{eq:det-ratio}
\det A/\det A_{\b,\b}= \det \,(a^\b_{kl})^{\phantom{\b}}_{k,l\in
\a}
\end{equation}
with $a^\b_{kl} = \det A_{\b\cup\{k\}, \b\cup\{l\}}/\det A_{\b,\b}$.
\end{lem}

\Proof If  $\G=\a\cup\b$ and $A_{\b,\b}$ is invertible we can write
\[
A=\begin{pmatrix}A_{\a,\a}&A_{\a,\b}\\A_{\b,\a}&A_{\b,\b}\end{pmatrix} =
\begin{pmatrix}A_{\a,\a}-A_{\a,\b}A_{\b,\b}^{-1}A_{\b,\a}&
A_{\a,\b}A_{\b,\b}^{-1}\\0&I\end{pmatrix}
\begin{pmatrix}I&0\\ A_{\b,\a}&A_{\b,\b}\end{pmatrix}\,,
\]
where $I$ is the identity matrix. Hence
\begin{equation}\label{eq:det1}
\det A /\det A_{\b,\b}=
\det(A_{\a,\a}-A_{\a,\b}A_{\b,\b}^{-1}A_{\b,\a})
\end{equation}
Now, if $A$ is Hermitian and positive definite then
$A_{\a,\b}A_{\b,\b}^{-1}A_{\b,\a}\ge0$, and the determinant is monotone
with respect to the
operator ordering. This implies \eqref{eq:det2} for invertible $A_{\b,\b}$. For general $A_{\b,\b}$ one can approximate
$A_{\b,\b}$ by invertible matrices. To prove \eqref{eq:det3} one can assume that
$\det A_{\b,\b}>0$ and then divide both sides
by $\det^2 A_{\b,\b}$. In view of \eqref{eq:det1}, \eqref{eq:det3} then takes the form
\[
\begin{split}
&\det  (A_{\a\cup\g,\a\cup\g}-A_{\a\cup\g,\b}A_{\b,\b}^{-1}A_{\b,\a\cup\g}) \\
&\qquad\qquad\le
\det(A_{\a,\a}-A_{\a,\b}A_{\b,\b}^{-1}A_{\b,\a})\,
\det(A_{\g,\g}-A_{\g,\b}A_{\b,\b}^{-1}A_{\b,\g})
\end{split}
\]
which follows from \eqref{eq:det2}. To verify \eqref{eq:det-ratio} we
have to show that the right-hand sides of \eqref{eq:det1} and
\eqref{eq:det-ratio} are identical. However, applying \eqref{eq:det1}
to the matrix $A_{\a\cup\{k\}, \a\cup\{l\}}$ it is easily seen that $a^\b_{kl}$
is precisely the $kl$-entry of the matrix
on the right-hand side of \eqref{eq:det1}.\EndProof

The second subject of this appendix is the proper choice of
integral kernels. Note that we need to evaluate such integral
kernels at the points of configurations in $\SP$. This leads us to
the following concept. Let us say that two functions
$t_1,t_2:\SP^2\to\C$ are \emph{indistinguishable} ($L$-almost
everywhere) if
\[
L\Big(\a\in\cX_0: t_1(x,y)\ne t_2(x,y)\text{ for some }x,y\in\a\Big)=0\,.
\]
The lemma below shows the significance of trace formulas for indistinguishability.
As before, we write $T_\L=P_\L TP_\L$.

\begin{lem}\label{lem:indistinguishable}
Let $T$ be a local trace class integral operator on $L^2(\SP,\l)$ and $t_1,t_2$
be two integral kernels of $T$ satisfying
\[
\Tr\, T_\L=\int_{\L}t_i(x,x)\,\l(dx)
\]
for all $\L\in\cB_0$, $i=1,2$. Then $t_1$ and $t_2$ are indistinguishable.
\end{lem}

\Proof Since both $t_1$ and $t_2$  are integral kernels of $T$, we have
$t_1=t_2$ almost everywhere with respect to $\l^{\otimes 2}$.
On the other hand, since the local trace formula holds for each $\L\in\cB_0$,
it follows that $t_1=t_2$ almost everywhere also with respect to
$\bar\l:= \l(x: (x,x)\in \,\cdot\,)$. So, the set $\cN=\{t_1\ne t_2\}$ is a
nullset for both $\l^{\otimes 2}$ and $\bar\l$.
This proves the lemma because
\[
\hat\cN:= \{\a\in\cX_0: (x,y)\in\cN\text{ for some }x,y\in\a\}
\]
is then a nullset for $L$.\EndProof 

As for the existence of suitable integral kernels, we consider
first the case of trace class operators. (For $\SP=\R^d$, the
trace formula \eqref{eq:trace-formula} appears also in \cite[Lemma
2]{So}.) Recall that a function $t:\SP\times\SP\to\C$ is called
\emph{positive definite} if for any $n\in \N$,
$\{c_i\}_{i=1}^n\subset\C$ and $\{x_i\}_{i=1}^n\subset \SP$,
\begin{equation}\label{eq:positivity}
\sum_{i,j=1}^n \ol{c}_i\,t(x_i,x_j)\,c_j\ge 0.
\end{equation}
\begin{lem}\label{lem:positive_kernel}
Let $T\ge 0$ be any trace class operator on $L^2(\SP,\l)$. Then
$T$ admits a positive definite integral kernel $t(x,y)$ satisfying
\begin{equation}\label{eq:trace-formula}
\Tr(P_\L T^k P_\L)=\int_{\SP^k}1_\L(x_1)\,t(x_1,x_2)\cdots
t(x_k,x_1)\,\l(dx_1)\cdots\l(dx_k)
\end{equation}
for all $k\ge 1$ and $\L\in\cB$.
\end{lem}

\Proof  Since $T$ is positive and of trace class, there exists a
Hilbert-Schmidt operator $S$ such that $T=S^*S$. Let $s(x,y)$ be
an integral kernel function of $S$ and $s^*(x,y):=\ol{s(y,x)}$ the
associated kernel of $S^*$. The Hilbert-Schmidt norm of $S$
satisfies $\|S\|_2^2=\int_{\SP^2}|s(x,y)|^2\l(dx)\l(dy)<\infty$.
There exists therefore a $\l$-nullset $\cN \subset\SP$ such that
$s(\cdot,y)\in L^2(\SP,\l)$ for $y\notin \cN $. Define $t(x,y)=\int
s^*(x,z)s(z,y)\,\l(dz)$ if $x,y\notin \cN $ and $t(x,y)=0$ otherwise.
Then $t$ is an integral kernel of $T$, and for
$\{c_i\}_{i=1}^n\subset\C$ and $\{x_i\}_{i=1}^n\subset \SP$ one
has
\[
\sum_{i,j=1}^n \ol{c}_i\, t(x_i,x_j)c_j=\int \big|\sum_{j=1}^n
c_j\,s(z,x_j)\big|^2\,\l(dz)\ge 0.
\]
On the other hand,
\[
\Tr\,T=\Tr\,S^*S=\|S\|_2^2=\int_{\SP^2}|s(x,y)|^2\l(dx)\l(dy)=\int_\SP
t(x,x)\l(dx),
\]
proving \eqref{eq:trace-formula} for $k=1$ and $\L=\SP$. To check
the case $\L\sub\SP$, it is sufficient to note that $SP_\L$ admits
the kernel $s(x,y)1_\L(y)$. The case of powers $T^k$ with $k>1$
follows similarly because $T^k=(S^{(k)})^*S^{(k)}$ with
$S^{(k)}=S^*S\cdots S^\sharp$, where $S^\sharp=S$ if $k$ is even
and $S^\sharp=S^*$ if $k$ is odd.\EndProof

Next we need to extend Lemma \ref{lem:positive_kernel} to integral
operators which are only local trace class.

\begin{lem}\label{lem:loc_trace_formula}
Let $T\ge 0$ be a bounded local trace class integral operator
on $L^2(\SP,\l)$. Then one can choose its integral
kernel $t(x,y)$ such that the following properties hold:
\begin{enumerate}
\item $t(\a,\a)\ge 0$ for $L$-a.a. $\a\in \cX_0$;
\item $t(x,y)$ is a Carleman
kernel, i.e., $t_x:=t(\cdot,x)\in L^2(\SP,\l)$ for $\l$-a.a. $x\in \SP$;
\item For any $\L\in \cB_0$, $\Tr(T_\L)=\int_\L t(x,x)\,\l(dx)$ and
\[
\Tr(P_\L T^k P_\L)=\int_\L (t_x,T^{k-2} t_x)\,\l(dx)
\]
for $k\ge 2$, where $(\cdot,\cdot)$ stands for the inner product in $L^2(\SP,\l)$.
\end{enumerate}
\end{lem}

\Proof
For each $\D\in\cB_0$ we choose an integral kernel $t_\D$ of
$T_\D$ according to Lemma \ref{lem:positive_kernel}. For any
$\L\sub\D$, $1_\L(x)t_\D(x,y)1_\L(y)$ is also an integral kernel
of $T_\L$ which satisfies the assumptions of Lemma
\ref{lem:indistinguishable}. As in the proof of that lemma, we
thus find a set $\cN_{\L,\D}$ which is a nullset for both
$\l^{\otimes 2}$ and $\bar\l$ such that $t_\L=t_\D$ on
$\L^2\sm\cN_{\L,\D}$.

Now let $(\D_n)_{n\ge1}$ be a sequence in $\cB_0$ increasing to
$\SP$, and consider the set $\cN=\bigcup_{n\ge1}\cN_{\D_n,\D_{n+1}}$.
For $(x,y)\notin\cN$ we let $t(x,y)$ be the common value of all
$t_{\D_n}(x,y)$ with $\D_n\ni x,y$, and define $t=0$ on $\cN$.
Since $T$ is supposed to admit an integral kernel $\tilde t$, it is then
clear that
\[
\int_{\SP^2}\ol{f(x)}\,\big[t(x,y)-\tilde t(x,y)\big] \,g(y)\, \l(dx)\l(dy) = 0
\]
whenever $f,g\in L^2(\SP,\l)$ vanish outside some $\D_n$, and
hence that $t=\tilde t$ $\l^{\otimes 2}$-a.e., meaning that $t$ is
also a kernel for $T$. By construction, $t$ is indistinguishable
from $t_{\D_n}$ under $L_{\D_n}$, and $t_{\D_n}$ is positive
definite. This yields assertion (a).

To prove (b) and (c) we proceed as follows.
Let $T_n:=T_{\D_n}$ for brevity. Since $T$ is bounded,
$P_\L T_n^k P_\L$ converges strongly to $P_\L T^k
P_\L$ as $n\to \infty$. In particular, it converges weakly.
Since also both $P_\L T_n^k P_\L $ and $P_\L T^k P_\L
$ are bounded from above by the trace class operator
$\|T\|^{k-1}T_\L$, it follows that the sequence $(P_\L
T_n^k P_\L)_{n\ge 1}$ converges to $P_\L T^k P_\L$
in trace norm \cite[Theorem 2.16]{Si}. Hence
\begin{equation}\label{eq:tr_convergence}
\Tr(P_\L T^kP_\L)=\lim_{n\to \infty}\Tr(P_\L T_n^kP_\L).
\end{equation}
However,  $t_{\D_n}$ has been chosen to satisfy the identity
\begin{equation}\label{eq:locTrace3}
 \Tr(P_\L T_n^kP_\L)=
\int_{\D_n^k}1_\L(x_1)\,t_{\D_n}(x_1,x_2)\cdots
t_{\D_n}(x_k,x_1)\,\l(dx_1)\cdots\l(dx_k),
\end{equation}
and we have seen that $t_{\D_n}$ can be replaced by $t$ 
in this integral. Hence, statement (c) will be proved
once we have shown that the right-hand side of 
\eqref{eq:locTrace3} converges to the integrals in (c).
For $k=1$, this is evident from the above. In the case 
$k=2$, \eqref{eq:locTrace3} means that
\[
\Tr(P_\L T_n^2
P_\L)=\int_{\SP^2}1_\L(x_1)1_{\D_n}(x_2)\,|t(x_1,x_2)|^2\,\l(dx_1)\l(dx_2).
\]
The integrand is nonnegative and increases to $1_\L(x_1)|t(x_1,x_2)|^2$ 
as $n\to \infty$. On the other hand,  $\Tr(P_\L T_n^2 P_\L)\le 
\|T\|\Tr(P_\L T_n P_\L)=\|T\|\Tr(T_\L)<\infty$ for all $n$ with
$\D_n\supset \L$. The monotone convergence theorem thus shows that $1_\L(x_1)|t(x_1,x_2)|^2$ is integrable on $\SP^2$. This proves statement (c)
for $k=2$, and also assertion (b) by Fubini's theorem.
Finally, let $k\ge3$. Eq.~\eqref{eq:locTrace3} can then be rewritten
in the form
\[
 \Tr(P_\L T_n^kP_\L)=
\int_{\L}(t_x,T_n^{k-2}t_x)\, \l(dx),
\]
and the integrand is bounded by $(t_x,T_n^{k-2}t_x)\le \|T\|^{k-2}(t_x,t_x)$
uniformly in $n$. We also know from the case $k=2$ that this 
upper bound is locally integrable.  We know further that 
$(t_x,T_n^{k-2}t_x)$ converges to
$(t_x,T^{k-2}t_x)$ for all $x$ with $t_x\in L^2(\SP,\l)$, and thus
for $\l$-a.a. $x$. Again, these functions are bounded by
$\|T\|^{k-2}(t_x,t_x)$. Together with \eqref{eq:tr_convergence}, 
this gives the trace formula in (c) for $k\ge3$. The proof is 
therefore complete.\EndProof

A final useful tool is the following projection-inversion lemma 
appearing without proof in \cite[p. 18]{OP}; see
\cite[Corollary 5.3]{ST2} for the matrix version.

\begin{lem}\label{lem:basic-order}
Let $T$ be any bounded positive operator with bounded inverse
$T^{-1}$. Then for any projection $P$,
\[
PT^{-1}P\ge P(PTP)^{-1}P.
\]
\end{lem}

\Proof Since $T$ is positive and invertible, any restriction of
$T$ to some subspace is invertible in this subspace. This means
that the operators $R:=P(PTP)^{-1}P$ and
$S:=P^\bot(P^\bot T^{-1}P^\bot)^{-1}P^\bot$ are
well-defined, where $P^{\bot}=I-P$. The key
observation is the decomposition formula
\begin{equation}\label{eq:decomposition}
PT^{-1}P=R+PT^{-1}ST^{-1}P\,.
\end{equation}
To see this we observe that $P^\bot T^{-1}(P^\bot+P)TP=0$.
Multiplying with $S$ from the left we find $P^\bot
TP=-S T^{-1}PTP$. Inserting this into the identity
$PT^{-1}(PTP+P^\bot TP)=P$ we get
\[
PT^{-1}\Big(I-ST^{-1}\Big)PTP=P.
\]
Multiplying with $R$ from the right and rearranging we
arrive at \eqref{eq:decomposition}. As the second operator on the
right-hand side of \eqref{eq:decomposition} is positive, the lemma
follows.\EndProof

\medskip\small
\noindent{\it Acknowledgments}. We are grateful to the
referees whose remarks helped to improve this paper.  H.J.Y. would
like to thank Prof. T. Shirai and Prof. Y. Takahashi for inviting
him to Kanazawa University. Part of this work was done during this
stay. We also thank Prof. Y. M. Park for informing us about
reference \cite{OP}. H.J.Y. was supported by Korea Research
Foundation Grant (KRF-2002-015-CP0038).



\begin{thebibliography}{99}
\bibitem{B} R. Bhatia, {\it Matrix analysis}, Springer, New York, 1997.
\bibitem{BOO} A. Borodin, A. Okounkov, and G. Olshanski,
Asymptotics of Plancherel measures for symmetric groups, {\it J.
Amer. Math. Soc.} {\bf 13} (3), 481-515 (2000).
\bibitem{BO} A. Borodin and G. Olshanski, Point processes and the
infinite symmetric group, Part VI: Summary of results, Available
at http://arxiv.org/abs/math.RT/9810015.
\bibitem{DV} D.J. Daley and D. Vere-Jones, {\it An introduction to the
theory of Poisson processes}, Springer-Verlag, New York, 1988.
\bibitem{GK} H.-O. Georgii and T. K\"{u}neth, Stochastic comparison of
point random fields, {\it J. Appl. Prob.} {\bf 34}, 868--881 (1997).
\bibitem{Gloetzl2} E. Gl\"otzl, Konstruktion der bedingten Energie eines
Punktprozesses, {\it Serdica} {\bf 7}, 217--233 (1981)
\bibitem{G} G. Grimmett, {\it Percolation}, 2nd ed., Springer, Berlin,
1999.
\bibitem{HS} P. R. Halmos and V. S. Sunder, {\it Bounded integral operators on $L^2$
spaces}, Springer-Verlag, Berlin etc, 1978.
\bibitem{KK} Y. G. Kondratiev and T. Kuna, Correlation functions
for Gibbs measures and Ruelle bounds, {\it Methods Funct. Anal.
Topology} {\bf 9}, 9-58 (2003).
\bibitem{Kozlov} O.K. Kozlov,  Gibbsian description of point random
fields, {\it Theory Prob. Appl.} {\bf 21}, 339--355 (1976.
\bibitem{L} R. Lyons, Determinantal probability measures,
{\it Publ. Math. Inst. Hautes \'Etudes Sci.} {\bf 98}, 167-212
(2003).
\bibitem{LS} R. Lyons and J. E. Steif, Stationary determinantal
process: Phase multiplicity, Bernoullicity, entropy, and
domination, {\it Duke Math. J.} {\bf 120} (3), 515-575
(2003).
\bibitem{M} O. Macchi, The coincidence approach to stochastic point
processes, {\it Adv. Appl. Prob.} {\bf 7}, 83--122 (1975).
\bibitem{MWM} K. Matthes, W. Warmuth, J. Mecke, Bemerkungen zu einer
Arbeit von Nguyen Xuan Xanh und Hans Zessin, {\it Math. Nachr.} {\bf
88}, 117--127 (1979).
\bibitem{Mecke} J.~Mecke, Station\"are zuf\"allige Ma{\ss}e auf
lokalkompakten abelschen Gruppen, {\it Z. Wahr\-schein\-lichkeitstheorie
verw. Geb.} {\bf 9}, 36--58 (1967).
\bibitem{MR} R. Meester and R. Roy, {\it Continuum percolation},
Cambridge University Press, 1996.
\bibitem{NZ} Nguyen X.X. and H. Zessin, Integral and differential
characterizations of the Gibbs process, {\it Math. Nachr.} {\bf 88},
105--115 (1979).
\bibitem{OP} M. Ohya and D. Petz, {\it Quantum entropy and its
use}, Springer-Verlag, Berlin, 1993.
\bibitem{Papangelou} F. Papangelou, The conditional intensity of general
point processes and an application to line processes, {\it Z.
Wahrscheinlichkeitstheorie verw. Geb.} {\bf 28}, 207--226 (1974)
\bibitem{Pr} C. J. Preston, Spatial birth-and-death processes, {\it
Bull. Inst. Int. Statist.} {\bf 46}, 371--391 (1976).
\bibitem{PrLNM} C. J. Preston, {\it Random Fields}, Lecture Notes in
Mathematics Vol. 534, Springer Verlag, Berlin etc., 1976.
\bibitem{RS} M. Reed and B. Simon, {\it Methods of modern mathematical
physics II: Fourier analysis, self-adjointness}, Academic Press,
New York etc, 1975.
\bibitem{Sh} H. Shimomura, Poisson measures on the configuration space
and unitary representations of the group of diffeomorphisms, {\it J.
Math. Kyoto Univ.} {\bf 34} (3), 599-614 (1994).
\bibitem{ST1}T. Shirai and Y. Takahashi, Random point field associated
with certain Fredholm determinant I: fermion, Poisson, and boson point
processes, {\it J. Funct. Anal.} {\bf 205}, 414--463 (2003).
\bibitem{ST2}T. Shirai and Y. Takahashi, Random point field associated
with certain Fredholm determinant II: fermion shift and its ergodic and
Gibbs properties, {\it Ann. Prob.} {\bf 31}, 1533--1564 (2003).
\bibitem{SY}T. Shirai and H. J. Yoo, Glauber dynamics for fermion point
processes, {\it Nagoya Math. J.} {\bf 168}, 139-166 (2002).
\bibitem{Si} B. Simon, {\it Trace ideals and their applications},
Cambridge University Press, Cambridge etc., 1979.
\bibitem{So} A. Soshnikov, Determinantal random point fields, {\it Russ.
Math. Surv.} {\bf 55}, 923-975 (2000).
\bibitem{Sp1} H. Spohn, Interacting Brownian particles: A study of
Dyson's model. In: G. Papanicolaou (ed.), {\it  Hydrodynamic behaviour
and interacting particle systems (Minneapolis, Minn., 1986)}, IMA Vol.
Math. Appl. 9 (1987).
\bibitem{Sp2} H. Spohn, Tracer dynamics in Dyson's model of interacting
Brownian particles, {\it J. Stat. Phys.} {\bf 47}, 669-679 (1987).
\bibitem{Strassen} V. Strassen, The existence of probability measures
with given marginals, {\it Ann. Math. Statist.} {\bf 36}, 423--439
(1965).
\bibitem{Y} H. J. Yoo, Gibbsianness of fermion random point fields,
Preprint (2003).
\end{thebibliography}
\end{document}